\documentclass{amsart}
\usepackage{amsmath,amsthm}
\usepackage{amsfonts,amssymb}

\newtheorem{theorem}{Theorem}[section]
\newtheorem{corollary}[theorem]{Corollary}

\newtheorem{proposition}[theorem]{Proposition}
\newtheorem{definition}[theorem]{Definition}

\theoremstyle{remark}

 \def\nb{{\mathbf n}}

 \def\CA{{\mathcal A}}
 \def\CD{{\mathcal D}}
 \def\CE{{\mathcal E}}

 \def\CL{{\mathcal L}}

 \def\CV{{\mathcal V}}

 \def\NN{{\mathbb N}}

 \def\RR{{\mathbb R}}
 
 \def\ZZ{{\mathbb Z}}

\newcommand{\wt}{\widetilde}

\input epsf.tex
 
\begin{document}

\title[difference equations and orthogonal polynomials]
{Second order difference equations and discrete orthogonal 
polynomials of two variables}  
\author{Yuan Xu}
\address{Department of Mathematics\\ University of Oregon\\
    Eugene, Oregon 97403-1222.}\email{yuan@math.uoregon.edu}

\date{\today}
\thanks{Work partially supported by the National Science Foundation under 
Grant DMS-0201669}
                         
\begin{abstract}
The second order partial difference equation of two variables 
\begin{align*} 
 & \CD u:= A_{1,1}(x) \Delta_1 \nabla_1 u + A_{1,2}(x) \Delta_1 \nabla_2 u 
     + A_{2,1}(x) \Delta_2 \nabla_1 u +  A_{2,2}(x) \Delta_2 \nabla_2 u\\
 & \qquad \qquad \qquad \qquad 
    +  B_1(x) \Delta_1 u  + B_2(x) \Delta_2 u = \lambda u,  \notag
\end{align*}
is studied to determine when it has orthogonal polynomials as solutions. 
We derive conditions on $\CD$ so that a weight function $W$ exists for
which $W \CD u$ is self-adjoint and the difference equation has polynomial 
solutions which are orthogonal with respect to $W$. The solutions are
essentially the classical discrete orthogonal polynomials of two variables.
\end{abstract}

\maketitle                      
 
\section{Introduction}
\setcounter{equation}{0}

Let $\Delta$ and $\nabla$ denote the forward and the backward difference 
operators, respectively, 
$$
\Delta f(t) = f(t+1) - f(t) \quad\hbox{and}\quad  \nabla f(t) = f(t) - f(t-1).
$$ 
The second order difference equation in one variable takes the form
\begin{equation} \label{eq:1.1}
  a(t) \Delta \nabla u+ b(t) \Delta u = \lambda u, \qquad \lambda \in \RR.  
\end{equation}
In the case that $a$ and $b$ are polynomials, it is known that some of 
the difference equations can have orthogonal polynomials as solutions. 
Such difference equations have been characterized completely and the 
orthogonal polynomials so derived are the classical discrete orthogonal 
polynomials (Hahn, Meixner, Krawtchouk and Charlier polynomials), see 
\cite{Al} for a survey in this direction and for references. 

The purpose of the present paper is to consider the second order 
difference equations in two variables. Throughout this paper we use the 
notation
$$
  e_1 = (1,0) \qquad \hbox{and} \qquad e_2 = (0,1).
$$ 
For $x = (x_1,x_2) \in \RR^2$ we define the forward and the backward 
partial difference operators by 
\begin{equation} \label{eq:1.2}
\Delta_i u(x) = f(x+e_i) - f(x) \quad\hbox{and}\quad 
\nabla_i u(x) = f(x) - f(x-e_i), \quad i = 1,2,
\end{equation}
respectively. We consider the second order partial difference equation  
\begin{align} \label{eq:1.3}
 & \CD u:= A_{1,1}(x) \Delta_1 \nabla_1 u + A_{1,2}(x) \Delta_1 \nabla_2 u 
     + A_{2,1}(x) \Delta_2 \nabla_1 u +  A_{2,2}(x) \Delta_2 \nabla_2 u\\
 & \qquad \qquad \qquad \qquad 
    +  B_1(x) \Delta_1 u  + B_2(x) \Delta_2 u = \lambda u,  \notag
\end{align}
where $A_{i,j}$ and $B_i$ are polynomials and $\lambda$ is a real
number. The goal is to characterize those equations that have orthogonal 
polynomials as solutions. 

The orthogonal polynomials that satisfy \eqref{eq:1.3} will be discrete
orthogonal polynomials of two variables. A sequence of polynomials, 
$\{p_\alpha\}$, of $d$ variables is called a sequence of discrete orthogonal
polynomials if there exist a set of isolated points $V \subset \RR^d$ and
a weight function $W$ defined on $V$ such that $p_\alpha$ are orthogonal 
with respect to the bilinear form
$$
  \langle f, g\rangle = \sum_{x \in V} f(x) g(x) W(x).
$$
Since we will discuss orthogonal polynomials associated with the difference 
equation \eqref{eq:1.3}, the set $V$ will be a lattice set, that is, a 
subset of $\ZZ \times \ZZ$. We should mention that a simple change of variables
$x \to h x$ shows that we can replace the difference operators $\Delta_i$ 
by the divided difference 
$$
[f(x_1+h,x_2) - f(x_1,x_2)]/{h} \quad \hbox{and} \quad 
[f(x_1,x_2+k) - f(x_1,x_2)]/{k}. 
$$
In this case the second difference $\Delta_i \nabla_j$ becomes 
second order divided difference; for example, $\Delta_1 \nabla_1$ and
$\Delta_1 \nabla_2$ become
\begin{align*}
& [f(x_1+h,x_2) - 2 f(x_1,x_2) + f(x_1 -h,x_2)]/{h^2} \quad \hbox{and} \\ 
& [f(x_1+h,x_2) - f(x_1+h,x_2-k) - f(x_1,x_2)+ f(x_1,x_2-k)]/{h k}, 
\end{align*}
respectively. In other words, the difference equations \eqref{eq:1.3}
is the model for the uniform grid. 

After the dilation indicated above, letting $(h, k) \to (0,0)$, the 
divided difference operators become the differential operators; in 
particular, $\Delta_i$ becomes $\partial/\partial_i$ and $\Delta_i \nabla_j$ 
becomes $\partial^2 /\partial_i \partial_j$. This way the difference 
equation becomes a differential equation. In this 
connection, it should be pointed out that the classification of the 
second order partial differential equations that have orthogonal
polynomial as solutions was carried out in \cite{KS}, see  
\cite{BSX,Kwon,Kwon2,Litt,Su} for further work in this direction. 

Comparing to the theory in one variable, the structure of discrete 
orthogonal polynomials in several variables is much more complicated.
Some basic results are obtained in \cite{X04}; the relevant ones will
be recalled in the following section. 

In the next section we will also discuss for what $A_{i,j}$ and $B_i$ the
equation \eqref{eq:1.3} has polynomial solutions. Such equations are
called admissible. In Section 3 we define and discuss the self-adjoint
form of the second order difference equations and give necessary
and sufficient conditions for the existence of a function $W$ so that
$W \CD$ is self-adjoint. We also define a notion of $W$ being consistent
with the difference equation, and show that consistence and self-adjointness
will imply that the difference equation has polynomial solutions that are
orthogonal with respect to $W$. In Section 4, we identify difference 
equations that have discrete orthogonal polynomials as solutions. Up to 
the mild restriction that we put on the difference equations, we believe
that these are the only such difference equations. The polynomials are 
the classical discrete orthogonal polynomials of two variables (see, 
for example, \cite{KM, Tr1,Tr2}) on the set $V = \{x: x_1 \ge 0, 
x_2 \ge 0, x_1+x_2 \le N\}$ and on the set $V = \{x: 0 \le x_1 \le M, 
0 \le x_2 \le N\}$, where both $N$ and $M$ can be infinity.  

\section{Preliminary and Admissible equations}
\setcounter{equation}{0}

\subsection{Discrete orthogonal polynomials}

Let $V$ be a subset of $\RR^d$. Orthogonal polynomials on $V$ depend 
on the structure of the polynomial ideal $I(V):= \{p \in \RR[x]: p(x) =0, 
\forall x \in V\}$ that has $V$ as its variety. The discrete orthogonal 
polynomials on $V$ can only consist of polynomials that do not belong to 
$I(V)$. Such polynomial subspaces have been identified and the basic 
structure of discrete orthogonal polynomials of several variables, including 
three-term relation and Favard's theorem, has been studied in \cite{X04}.

Let $\RR[V]$ denote the subspace of polynomials. This is the space that
orthogonal polynomials belong to. Then $\RR[V] \cong 
\RR[x_1,\ldots, x_d]/I(V)$.  There is a lattice set $\Lambda = \Lambda(V)$
such that every polynomial $P \in \RR[V]$ can be written as
$$
  P(x) = \sum_{\alpha \in \Lambda} c_\alpha x^\alpha  \mod I(V), 
           \qquad c_\alpha \in \RR.
$$  
One particular result in \cite{X04} shows that the set $\Lambda$ 
satisfies the following property
\begin{equation*} 
\alpha \in \Lambda \quad \hbox{implies} \quad \alpha-\beta\in \Lambda, \quad 
\hbox{whenever  $\alpha - \beta \in \NN_0^d$ and $\beta \in \NN_0^d$}.
\end{equation*}
For two variables, $d =2$, this means that $\Lambda$ must be of a stair
shape; more precisely, there is a sequence of positive integers $n_i$,
which satisfies $n_m \le n_{m-1}\le \ldots\le n_0$ (some of the $n_i$ 
can be positive infinity and so can $m$), such that 
\begin{equation} \label{eq:2.1}
 \Lambda = \{(k,l): 0 \le l \le m, 0 \le k \le n_l\}.  
\end{equation} 
Hence, any polynomial $P \in \RR[V]$ can be written as 
$$
  P(x) = \sum_{l=0}^m \sum_{j=0}^{n_l} c_{l,j} x_1^l x_2^j, \qquad
    c_{l.j} \in \RR. 
$$
We will assume that $\Lambda$ contains $\{(i,j): 0 \le i+j \le 2 \}$ and 
at least one of the points $(2,1)$ and $(1,2)$. 

Given $\Lambda$ as in \eqref{eq:2.1}, the highest degree of monomials 
appeared in $\RR[V]$ is denoted by $\nb$, which could be positive infinity.
It follows that  
\begin{equation*}  
  \nb: = \max \{ n_l + l : 0 \le l \le m \}.
\end{equation*} 
For $\Lambda = \Lambda(V)$, denote $\Lambda_k(V) = \{(i,j): i+j =k\}$ for 
each $k < \nb$. Let $r_k$ denote the distinct elements in $\Lambda_k(V)$, 
which is the number of monomials of degree exactly $k$ in $\RR[V]$. For 
each $j$ satisfying $0 \le j \le \nb$, it can be verified that 
\begin{equation} \label{eq:2.2}
  r_k = k + 1 - \sum_{l = 0}^m (k - l - n_l)_+ , 
\end{equation} 
where $(a)_+ = 0$ if $a \le  0$ and $(a)_+ = a$ if $a > 0$. 

Let $W$ be a real function on $V$ and $W(x) \ne 0$ on $V$ for any $x \in V$.
Assume that 
$$
  \sum_{x \in V} |x_1^i x_2^j| |W(x)| < \infty \qquad \hbox{for all 
    $(i,j) \in \NN_0^2$}
$$ 
in the case $V$ is an infinite set. Let $\Pi^2 = \RR[x_1,x_2]$ and let 
$\Pi_n^2$ denote the subspace of polynomials of degree at most $n$. 
Define the bilinear form $\langle \cdot, \cdot\rangle$ on 
$\Pi^2 \times \Pi^2$ by 
$$
\langle f,g \rangle = \CL (f g), \quad \hbox{where} \quad
 \CL(f):= \sum_{x \in V} f(x)W(x). 
$$
If $\langle f,g \rangle =0$, we say that $f$ and $g$ are orthogonal to each 
other with respect to $W$ on the discrete set $V$. A polynomial of degree
$k$ is said to be an orthogonal polynomial with respect to $W$ if it is 
orthogonal to all polynomials of degree lower than $k$ in $\RR[V]$. 
For a given weight function $W$, let $\CV_k$ denote the space of orthogonal
polynomials of degree exactly $k$ on $V$. The dimension of $\CV_k$ is $r_k$ 
in \eqref{eq:2.2}. Only in the case of $V=\{x: x_1 \ge 0,x_2 \ge 0, 
0 \le x_1+x_2 \le N\}$ this dimension is $r_k = k+1$, which is equal to 
the dimension for the space of continuous orthogonal polynomials 
(see \cite{DX}).

Since we are interested in discrete orthogonal polynomials that are solutions
of the difference equation \eqref{eq:1.3}, the set $V$ is a lattice set. 
If $u$ satisfies the difference equation \eqref{eq:1.3} then the function
$u^* = u(x_1-k,x_2-l)$ will satisfy a similar difference equation. Hence We 
can assume that $(0,0)$ belong to $V$. For a given $(i,j) \in \ZZ\times \ZZ$, 
let $B_{i,j}$ denote the square that has $(0,0), (i,0), (0,j)$ and $(i,j)$ as 
corners. In order that the differences make sense, we further assume that $V$
satisfies the following property:  
\begin{equation} \label{eq:2.3} 
  (i,j) \in V \quad \hbox{implies} \quad (k,l) \in V, \quad 
  \hbox{whenever $(k,l) \in B_{i,j}$}.
\end{equation}
As pointed out in \cite{X04}, in some cases, the index set $\Lambda$ can 
be the same as the set $V$. This happens, in particular, in the following
two important cases, 

\par\noindent
{\bf Type A}. $V = \{(x_1,x_2): 0 \le x_1 \le M, 0 \le x_2 \le N\}$,

\par\noindent
{\bf Type B}. $V = \{(x_1,x_2): x_1 \ge 0, x_2 \ge 0, 0 \le x_1 +x_2 \le N\}$,

\noindent 
that are of particular interest to us.  We note that both $N$ and $M$ 
can be infinity. 

\subsection{Admissible difference equations}

Since the forward and the backward difference equations are related by 
\begin{equation} \label{eq:2.4}
  \Delta_j  = \nabla_j \Delta_j + \nabla_j, 
\end{equation}
a second order partial difference operator can be expressed in different
forms. We shall fix one particular form and call it the standard form.

\begin{definition}
Let $A_{i,j}$ and $B_i$, $i,j = 1,2$, be polynomials of two variables. 
We call the operator $\CD$ defined by 
\begin{align} \label{eq:2.5}
\CD: = & A_{1,1}(x) \Delta_1 \nabla_1 + A_{1,2}(x) \Delta_1 \nabla_2 
         + A_{2,1}(x) \Delta_2 \nabla_1  +  A_{2,2}(x) \Delta_2 \nabla_2\\
      & +  B_1(x) \Delta_1  + B_2(x) \Delta_2  \notag
\end{align}
a standard second order partial difference operator. 
\end{definition}

An alternative way of writing a second order partial difference operator is 
\begin{align*} 
\wt \CD:= & \wt A_{1,1}(x)\Delta_1 \nabla_1 + \wt A_{1,2}(x)\Delta_1 \nabla_2 
  + \wt A_{2,1}(x) \Delta_2 \nabla_1  + \wt A_{2,2}(x) \Delta_2 \nabla_2\\
 & + \wt B_1(x) \nabla_1 + \wt B_2(x) \nabla_2,  \notag
\end{align*}
which can be easily converted to the standard form by using \eqref{eq:2.4}.

\begin{definition}
The equation $\CD u = \lambda u$ is called admissible on $V$ if for any 
$k\in \NN_0$, there is a number $\lambda_k$ such that the equation 
$\CD u = \lambda_k u$ has $r_k$ linearly independent polynomial solutions 
and it has no non-trivial solutions in the set of polynomials of degree 
small than $k$.
\end{definition}

The following proposition shows that we can assume that the linearly
independent polynomial solutions are given by those polynomials whose 
highest term contains a single monomial. 

\begin{proposition}
The equation $\CD u = \lambda u$ is admissible in $\RR[V]$ if and only if,
for every integer $k$ that satisfies $0 \le k \le n$, there exists a number 
$\lambda_k$ such that the equation $\CD u = \lambda_k u$ has $r_k$ linearly 
independent solutions in the form of 
\begin{equation} \label{eq:2.6}
  P_{k,l}(x) = x_1^k x_2^l + R_{k,l}(x), \qquad (k,l) \in \Lambda_k,
\end{equation}  
where $R_{k,l} \in \Pi_k^2$, and the equation has no non-trivial solutions
in $\RR[\Lambda]$ of degree less than $k$.
\end{proposition}

\begin{proof}
That the equation is admissible implies that, for each $k$, there exist
solutions $Q_{i,j} \in \RR[V]$, $\deg Q_{i,j} = k$. Introducing the 
notation 
$$
  x_\Lambda^k = \{x_1^i x_2^j: (i, j) \in \Lambda, \quad i+j \le k\}
$$
and regard it as a column vector, we can write 
$$
  Q_{i,j} (x) = a_{i,j}^T x_\Lambda^k + R_{i,j}, \qquad \deg R_{i,j} < k 
$$
where $a_{i,j}^T$ is a row vector of real numbers. We show that 
$\{a_{i,j}^T x_\Lambda^k: (i,j) \in \Lambda \}$ is linearly independent. 
Indeed, suppose otherwise, then there exist $c_{i,j}$ such that
$\sum_{(i,j) \in \Lambda} c_{i,j} a_{i,j}^T x_\Lambda^k =0$, which 
implies that 
$$
 \sum_{(i,j) \in \Lambda} c_{i,j} Q_{i,j} =  
  \sum_{(i,j) \in \Lambda} c_{i,j} R_{i,j} : = R_k \in \Pi_{k-1}^2.
$$ 
Since $Q_{i,j}$ satisfies $\CD u = \lambda_k u$, we conclude that
$$
 \CD R_k = \CD \sum_{(i,j) \in \Lambda} c_{i,j} Q_{i,j} = 
   \lambda_k  \sum_{(i,j) \in \Lambda} c_{i,j} Q_{i,j} = \lambda_k R_k
$$
so that $R_k \equiv 0$ by the definition of admissibility. Consequently,
since the set of $\{Q_{i,j}\}$ is linearly independent, $c_{i,j} = 0$ for
$(i,j)\in \Lambda$. Thus, $\{a_{i,j}^T x_\Lambda^k:(i,j)\in\Lambda_k\}$ is 
linearly independent and the matrix $(a_{i,j}^T)_{(i,j)\in \Lambda_k}$ is
invertible. Let 
$$
(P_{i,j})_{(i,j)\in \Lambda_k}:= 
 \left[(a_{i,j}^T)_{(i,j)\in \Lambda_k}\right]^{-1}
  (Q_{i,j})_{(i,j)\in \Lambda_k};
$$
then $P_{i,j}$ are polynomials of the form $P_{i,j}(x)=x_1^i x_2^j+ p_{i,j}$
where $\deg p_{i,j} < i+j = k$. 
This proves one direction of the proposition. The other direction is trivial,
since the set $\{P_{i,j}\}$ in \eqref{eq:2.6} is evidently linear independent. 
\end{proof}

\begin{proposition}
For the equation $\CD u = \lambda u$ to be admissible it is necessary 
that $A_{i,j}$ are polynomials of degree at most $2$ and $B_i$ are 
polynomials of degree at most $1$.
\end{proposition}

\begin{proof}
Suppose that $\CD u = \lambda u$ is admissible. Then there is a system of
polynomials, $P\{_{i,j}\}$ as in \eqref{eq:2.6}, which are solutions
of the equation. Let 
$A_{i,j}^0$ and $B_i^0$ denote the part of $A_{i,j}$ and $B_i$ that has 
degree larger than $2$, respectively. We need to prove that $A_{i,j}^0(x)
\equiv 0$ and $B_i^0(x) \equiv 0$. Substituting $P_{i,j}$ into 
$\CD u = \lambda u$ and examing the highest terms of the resulting equation, 
we obtain
\begin{align*}
 & A_{1,1}^0(x) \Delta_1 \nabla_1 x_1^i x_2^j + A_{1,2}^0(x) \Delta_1 \nabla_2 
   x_1^i x_2^j + A_{2,1}^0(x) \Delta_2 \nabla_1 x_1^i x_2^j \\
 & \qquad \qquad +  A_{2,2}^0(x) \Delta_2 \nabla_2 x_1^i x_2^j 
+  B_1^0(x) \Delta_1 x_1^i x_2^j + B_2^0(x) \Delta_2 x_1^i x_2^j=0  
\end{align*}
for $(i,j)\in \Lambda$. Let $(i,j) = (1,0)$ and $(i,j) = (0,1)$; then 
$\Delta_i \nabla_j x_1^i x_2^j =0$ and we conclude that 
$B_1^0(x) = B_2(x)^0 \equiv 0$. Next we let $(i,j) = (2,0)$; then only
$\Delta_1 \nabla_1 x_1^i x_2^j$ is nonzero and, consequently, 
$A_{1,1}^0(x) \equiv 0$. Similarly, we deduce $A_{2,2}^0(x) \equiv 0$ by 
setting $(i,j) = (0,2)$. Furthermore, the choice $(i,j) = (1,1)$ leads to 
$$
  A_{1,2}^0(x) + A_{2,1}^0(x) \equiv 0. 
$$
Finally, let $(i,j) = (1,2)$ (or $(2,1)$); we conclude that 
$$
  (2x_1+1) A_{1,2}^0(x) + (2x_1 -1) A_{2,1}^0(x) \equiv 0.
$$
Together, the last two displayed equations imply $A_{1,2}^0(x)=A_{2,1}^0(x) 
\equiv 0$. 
\end{proof} 

Since $A_{i,j}$ are polynomials of second degree and $B_i$ are polynomials
of the first degree, we can assume that they have the form. 
\begin{align*}
A_{i,j}(x)& = a_{i,j}^{2,0} x_1^2 + a_{i,j}^{1,1} x_1 x_2 +a_{i,j}^{0,2} x_2^2 
 + a_{i,j}^{1,0} x_1 + a_{i,j}^{0,1} x_2 +  a_{i,j}^{0,0}; \\
B_i(x) & = b_i^1 x_1 + b_i^2 x_2 + b_i^0, \qquad i,j = 1,2.   
\end{align*}

\begin{theorem} \label{thm:2.5}
Assume that the quadratic parts of $A_{1,2}$ and $A_{2,1}$ are equal; 
that is, $A_{1,2} - A_{2,1} \in \Pi_1^2$. 
Then the equation $\CD u = \lambda u$ is admissible if and only if 
\begin{equation}\label{eq:2.7}
 A_{i,j}(x) = a x_i x_j + a_{i,j} x_1 + b_{i,j} x_2 + c_{i,j}
\quad\hbox{and}\quad B_i = b x_i + d_i
\end{equation}
for $i,j = 1, 2$ with $c_{1,1}+c_{1,2}+c_{2,1}+ c_{2,2} =0$ and 
\begin{equation}\label{eq:2.8}
\lambda_k = k(k a - a +b), \qquad 0 \le k \le \nb,   
\end{equation} 
where $a$ and $b$ are real numbers such that 
$a p +  b \ne 0$ for all nonnegative integer $p$ satisfying $0 \le p \le \nb$.
\end{theorem}

\begin{proof}
Again consider the highest terms of the polynomials in $\CD P_{k,l} - 
\lambda P_{k,l}$. Since 
$\Delta_i x_i^j = (x_i + 1)^j - x_i^j = j x_i^{j-1} + \ldots$ and 
$\nabla_i x_i^j = x_i^j - (x_i - 1)^j  =  j x_i^{j-1} + \ldots$,
we obtain
\begin{align*}
& A_{1,1}(x) k(k-1) x_1^{k-2} x_2^l + A_{1,2}(x) k l
  x_1^{k-1} x_2^{l-1} + A_{2,1}(x) k l x_1^{k-1} x_2^{l-1} \\
&  \qquad \qquad + l(l-1)  A_{2,2}(x) x_1^k x_2^{l-2} 
  +B_1(x) k x_1^{k-1} x_2^l + B_2(x) l x_1^k x_2^{l-1} 
= \lambda_{j+k} x_1^k x_2^l.
\end{align*}
Using the formula of $A_{i,j}$ and $B_i$, we can then derive equations on the
coefficients of $A_{i,j}$ and $B_i$ by comparing the coefficients of 
$x_1^k x_2^l$. 

The coefficients of $x_1^{k-2}x_2^{l+2}$ and  $x_1^{k+2}x_2^{l-2}$ lead to 
equations 
$$  
   k(k-1) a_{1,1}^{0,2}=0 \quad \hbox{and} \quad   l(l-1) a_{2,2}^{2,0}=0
$$
which imply that $a_{1,1}^{0,2} = a_{2,2}^{2,0} =0$. Using the assumption 
that $a_{1,2}^{k,l}= a_{2,1}^{k,l}$ for $k + l =2$, the coefficients for
$x_1^{k-1} x_2^{l+1}$, $x_1^k x_2^l$ and $x_1^{k+1} x_2^{l-1}$ lead to the 
equations
\begin{align} \label{eq:2.9}
\begin{split}
k(k-1)a_{1,1}^{1,1}+ 2 k la_{1,2}^{0,2}+ k b_1^2 &=0,\\ 
k(k-1)a_{1,1}^{2,0}+ 2 k l a_{1,2}^{1,1}+ l(l-1) a_{2,1}^{1,1} + k b_1^1+ 
l b_2^2 &= \lambda_{k+l},\\
2 k l a_{1,2}^{2,0}+ l(l-1) a_{2,2}^{1,1} + l b_2^1 & = 0
\end{split}
\end{align}
for $(k,l) \in \Lambda$. In particular, if $(k,l) = (1,0)$ and $(k,l) = (0,1)$
then the first and the third equations of \eqref{eq:2.9} show that $b_1^2 =0$ 
and $b_2^1=0$, respectively. Furthermore, $k =1$ in the first equation   
shows that $a_{1,2}^{0,2} = 0$ and, consequently, $a_{1,1}^{1,1} = 0$. 
Similarly, $l =1$ in the third equation shows that $a_{1,2}^{2,0} =0$ and 
$a_{2,2}^{1,1}=0$. Setting $k=0$ or $l =0$ in the second equation of 
\eqref{eq:2.9} leads to the equations
$$
 l(l-1) a_{2,2}^{0,2} + l b_2^2  = \lambda_l \quad \hbox{and}\quad
 k(k-1) a_{1,1}^{2,0} + k b_1^1  = \lambda_k. 
$$
In particular, $k = l =1$ gives $b_2^2=b_1^1$ and, consequently,
$a_{1,1}^{2,0} = a_{2,2}^{0,2}$. Furthermore, letting $(k,l) = (1,1)$ and 
$(k,l) = (2,0)$ in the second equation of \eqref{eq:2.9} gives the equation
$$
  2 a_{1,2}^{1,1} + 2 b = \lambda_2 \quad \hbox{and} \quad 
  2 a_{1,1}^{2,0} + 2 b = \lambda_2 
$$
where we write $b = b_1^2 = b_2^1$. Consequently, we conclude that 
$a_{1,1}^{2,0} = a_{1,2}^{1,1}$. Writing $a = a_{1,1}^{2,0} = a_{1,2}^{1,1}
 = a_{2,2}^{0,2}$, we have proved that $A_{i,j}$ and $B_i$ are in the 
forms specified. Moreover, the second equation of \eqref{eq:2.9} with $l=0$
or with $k =0$ shows that
$$
 k(k-1) a +  k b = \lambda_k \quad\hbox{and}\quad  l(l-1)a + l b = \lambda_l.
$$
This is \eqref{eq:2.8} and that $\lambda_l \ne 0$ shows that $a p + b \ne 0$
for all nonnegative integer $p$ that satisfies $0 \le p \le \nb-1$.

On the other hand, assume that $A_{i,j}$ and $B_i$ are given as in 
\eqref{eq:2.7} and $\lambda_k$ is given as in \eqref{eq:2.8}. Let 
$(t)_k = t(t+1) \ldots (t+k-1)$ be the Pochhammer symbol. Define
$$
  m_k(t) = (t)_k /k! \quad\hbox{and}\quad m_{k,l}(x) = m_k(x_1) m_l(x_2). 
$$
For each pair of nonnegative integers $(r,s)$ satisfying $r+s =n$, 
substituting the polynomial 
$$
P_{r,s} (x) = m_{r,s}(x)  + \sum_{k+l < n} f_{k,l} m_{k,l}(x)
$$ 
into the difference equation $\CD u = \lambda_n u$ leads to a linear 
system of equations for the coefficients $f_{i,j}$. Indeed, the definition
of $m_{k,l}$ implies that
$$
 \Delta m_k(x) = m_{k-1}(x+1) = \sum_{j=0}^{k-1} m_{j,l}(x) 
   \quad \hbox{and} \quad \nabla m_{k}(x) = m_{k-1}(x), 
$$
from which it follows that
$$
 \Delta_1 \nabla_1 m_{k,l} = \sum_{j=0}^{k-2} m_{j,l}, 
   \quad \hbox{and} \quad   
 \Delta_1 \nabla_2 m_{k,l} = \sum_{j=0}^{k-1} m_{j,l-1}, 
$$
and similar formula for $\nabla_1 \Delta_2 m_{k,l}$ and
$\Delta_2 \nabla_2 m_{k,l}$. Furthermore, there is also $t m_k(t) =
(k+1) m_{k+1}(t) - k m_k(t)$, from which follows formulas such as 
\begin{align*}
x_1^2 \Delta_1 \nabla_1 m_{k,l}(x) & = k(k-1) m_{k,l} - (k-1)^2 m_{k-1,l}\\
x_1 x_2 \Delta_1 \nabla_1 m_{k,l}(x) & = k l m_{k,l}(x) - k(l-1) m_{k,l-1}(x).
\end{align*}
Using these formulas, a tedious computation shows that the coefficients of
$m_{r,s}$ in the expression $\CD P_{r,s} - \lambda_n P_{r,s}$ are
linear in $f_{k,l}$ with $r+s \le  k+l \le m$. Since the set
of $m_{k,l}$ is clearly linearly independent, these coefficients must 
be zero, which gives a linear system of equations in variables 
$f_{k,l}$. 

The coefficient of $m_{r,s}$ turns out to be $c_{1,1}+c_{1,2}+
c_{2,1}+c_{2,2}$, which is zero by assumption. Clearly, this also shows that
this condition is necessary. The coefficients of $m_{r-1,s}$ and $m_{r,s-1}$
lead to the following two equations
\begin{align*}
& (r-1)( -a (r-1) + a_{1,1})+s(-a (r-1)+a_{2,1})+s b_{1,2}+d_1 + 
   \lambda_{r-1,s} f_{r-1,s} =0, \\
& r(-a (s-1)+a_{1,2})+ (s-1)( -a (s-1) + a_{2,2}) +  
 r b_{2,1}+d_2 +  \lambda_{r,s-1} f_{r,s-1} =0, 
\end{align*}
from which $f_{r-1,s}$ and $f_{r,s-1}$ are determined uniquely. The 
coefficient of $m_{k,l}$ for $k+l < r+s$ is a linear combination of 
$f_{p,q}$ in which $p+q \ge k+l$ and $f_{k,l}$ is the only element with 
$p+q = k+l$. The coefficient of $f_{k,l}$ is 
$$
 \lambda_{k,l} - \lambda_{r,s} = (k + l - r - s) (b + a (-1 + k + l + r + s))
$$
which is nonzero by assumption. Hence, starting from $f_{r,s-1}$ and 
$f_{r,s-1}$ and works our way down, we conclude that $f_{k,l}$ is uniquely 
determined. 
\end{proof}

If $u$ satisfies an admissible difference equation $\CD u$, then under
a change of variables $x \mapsto x +h$, where $h \in \ZZ^2$, the function
$\wt u(x) = u(x + h)$ will satisfies a similar admissible difference equation.
We consider these equations equivalent. In particular, we can choose $h$ such
that $B_1(x) = b x_1$ and $B_2(x) = b x_2$. 

We note that an affine transform $x_1 = p_{1,1} x_1' + p_{1,2} x_2'$, 
$x_2 = p_{2,1} x_1' + p_{2,2} x_2'$, however, will alter the difference 
equations significantly.   

The characterization of partial differential equations in two variables 
that have orthogonal polynomials as solutions is given up to affine 
transformations; in other words, if $u(x)$ satisfies a differential 
equation, then the function $u^*(x') := u(x)$ under the affine transformation 
also satisfies a similar differential equation and the two equations are 
considered to be the same type. This allows us to use affine transformation 
to reduce each type of equations to its simplest form and proceeds from 
there. For the difference equations, however, the affine transformation
can no longer be used as a tool to reduce the type of equations. 

\section{Self-adjoint difference equation and the weight function}
\setcounter{equation}{0}

\subsection{Self-adjoint difference equations}
For the difference equation in one variable, the self-adjoint form is
$$
  \Delta ( \sigma(x) \nabla u(x)) - \lambda u(x) = 0.
$$

For difference equations in two variables, we introduce the following 
definition: 

\begin{definition} \label{defn:3.1}
Let $\sigma_{i,j}$ be continuous functions of two variables. The operator
$$
\CE u =
  \nabla_1 \left[ \sigma_{1,1}(x) \Delta_1 u +
  \sigma_{1,2}(x) \Delta_2 u \right] +  \nabla_2 
\left[ \sigma_{2,1}(x) \Delta_1 u + \sigma_{2,2}(x) \Delta_2 u \right] 
$$
is called a self-adjoint second order difference operator. The
equation $\CE u = \lambda u$ is called a self-adjoint difference equation.
\end{definition} 

\begin{proposition}\label{prop:3.3}
Let $\CD$ be the difference operator in \eqref{eq:2.5}. There is a function
$W$ such that $W(x) \CD$ is self-adjoint if and only if $W$ satisfies
\begin{align}\label{eq:3.1}
\begin{split}
 W B_1 & =\Delta_1(W A_{1,1})+\Delta_2 (W A_{1,2}) \\
 W B_2 & =\Delta_1(W A_{2,1})+\Delta_2 (W A_{2,2}). 
\end{split}
\end{align}
\end{proposition}

\begin{proof}
First we rewriting the self-adjoint form of $\CE u$ in the standard form of 
the difference operator $\CD u$. For this purpose the following formulas
are useful:
\begin{equation}\label{eq:3.2}
\Delta_i ( f g ) = f \Delta_i g + g(\cdot + e_i) \Delta_i f 
\quad\hbox{and}\quad
\nabla_i ( f g ) = f \nabla_i g + g(\cdot - e_i) \nabla_i f 
\end{equation}
where $i = 1,2$ and  
\begin{equation}\label{eq:3.3}
 \nabla_i u(\cdot+e_j) = \Delta_j \nabla_i u + \nabla_i u =
    \Delta_j \nabla_i u + \Delta_i u - \Delta_i \nabla_i u, \quad i,j = 1,2;
\end{equation} 
all can be easily verified. Using these formulas we can write 
\begin{align*}
 \CE u = & \sigma_{1,1}(\cdot - e_1) \Delta_1\nabla_1 u 
         + \sigma_{1,2}(\cdot - e_2)\Delta_1\nabla_2 u \\
       & + \sigma_{2,1}(\cdot - e_1)\Delta_2 \nabla_1 u 
         + \sigma_{2,2}(\cdot - e_2)\Delta_2\nabla_2 u\\
         & + (\nabla_1 \sigma_{1,1}+ \nabla_2 \sigma_{1,2}) \Delta_1 u 
           + (\nabla_1 \sigma_{2,1}+ \nabla_2 \sigma_{2,2}) \Delta_2 u.
\end{align*}
Comparing with the standard form $\CD u$ in \eqref{eq:2.5} and solving for 
$\sigma_{i,j}$ leads to 
\begin{align*}
& A_{1,1}=\sigma_{1,1}(\cdot-e_1), \quad A_{1,2} = \sigma_{1,2}(\cdot - e_2),
 \quad B_1 = \nabla \sigma_{1,1} + \nabla_2 \sigma_{1,2}, \\
& A_{2,1}=\sigma_{2,1}(\cdot - e_1), \quad A_{2,2} = \sigma_{2,2}(\cdot - e_2),
\quad B_2 = \nabla_1 \sigma_{1,2} +  \nabla_2 \sigma_{2,2}.
\end{align*}
Consequently, it follows that 
$$
   \sigma_{i,j} = A_{i,j}(\cdot + e_j), \qquad 1 \le i,j \le 2.
$$
Using the fact that $\nabla_i u(\cdot + e_i) = \Delta_i u$, we then obtain
$$
B_1 = \Delta_1 A_{1,1} + \Delta_2 A_{1,2} \quad\hbox{and}\quad 
B_2 = \Delta_1 A_{2,1} + \Delta_2 A_{2,2}.
$$
That is, the operator $\CD$ is self-adjoint only if the above relations
hold. Evidently, if $B_1$ and $B_2$ are given by the above relation, then
$\CD u$ can be written in the self-adjoint form. Hence, $\CD$ is self
adjoint if and only if these relations hold. Consequently, $W \CD$ 
is self-adjoint if and only if \eqref{eq:3.1} holds.
\end{proof} 

Using \eqref{eq:3.2}, the system of equations \eqref{eq:3.1} can be written as 
\begin{align}\label{eq:3.4}
\begin{split}
  A_{1,1}(x+e_1) W(x+e_1) + A_{1,2}(x+e_2) W(x+e_2) & = r_1(x) \\
  A_{2,1}(x+e_1) W(x+e_1) + A_{2,2}(x+e_2) W(x+e_2) & = r_2(x) 
\end{split}
\end{align}
where $r_i$ are polynomials defined by 
\begin{align*}
 r_1 = A_{1,1}+ A_{1,2} + B_1 \quad\hbox{and}\quad
 r_2 = A_{2,1}+ A_{2,2} + B_2.
\end{align*}
If the determinant of this system of equations is nonzero, then we can 
solve for $\nabla_1 W$ and $\nabla_2 W$. Let us define 
\begin{align*}
\alpha(x):= & A_{1,1}(x+e_1) A_{2,2}(x+e_2) - A_{1,2}(x+e_2)A_{2,1}(x+e_1),\\ 
\beta_1(x):= & r_1(x) A_{2,2}(x+e_2) - r_2(x) A_{1,2}(x+e_2), \\
\beta_2(x):= & r_2(x) A_{1,1}(x+e_1) - r_1(x) A_{2,1}(x+e_1).  
\end{align*} 

\begin{theorem}\label{thm:3.3}
Assume that $\alpha(x)\ne 0$ for $x \in V$. Then there is a $W$ such that 
$W \CD u$ is self-adjoint if and only if 
\begin{equation}\label{eq:3.5}
\beta_1 (x) \alpha(x+e_2)\beta_2 (x+e_1) = 
\beta_2 (x) \alpha(x+e_1) \beta_1 (x+e_2), \qquad x\in V. 
\end{equation}
\end{theorem}

\begin{proof}
Since $\alpha (x) \ne 0$ means the determinant of the system of equations
\eqref{eq:3.4} is nonzero, we can solve \eqref{eq:3.4} to get
\begin{equation}\label{eq:3.5a}
  \frac{W(x+e_1)}{W(x)} = \frac{\beta_1(x)}{\alpha(x)}  \quad\hbox{and}\quad
  \frac{W(x+e_2)}{W(x)} = \frac{\beta_2(x)}{\alpha(x)}.   
\end{equation}
These equations imply 
\begin{align*}
 W(x+e_1+e_2) & =  W(x + e_1) \frac{\beta_2(x+e_1)}{\alpha(x+e_1)} = 
    W(x)\frac{\beta_1(x)}{\alpha(x)} \frac{\beta_2(x+e_1)}{\alpha(x+e_1)},\\ 
 W(x+e_1+e_2) & =  W(x + e_2) \frac{\beta_1(x+e_2)}{\alpha(x+e_2)} = 
    W(x)\frac{\beta_2(x)}{\alpha(x)} \frac{\beta_1(x+e_2)}{\alpha(x+e_2)}.    
\end{align*}
Consequently, we conclude that 
$$
 \frac{\beta_1(x)}{\alpha(x)} \frac{\beta_2(x+e_1)}{\alpha(x+e_1)} =
 \frac{\beta_2(x)}{\alpha(x)} \frac{\beta_1(x+e_2)}{\alpha(x+e_2)}    
$$
from which the equation \eqref{eq:3.5} follows. 

On the other hand, assume that \eqref{eq:3.5} holds. We need to show that 
the system of equations \eqref{eq:3.4} has a solution.
We denote $\gamma_j:= \beta_j/\alpha$. The condition \eqref{eq:3.5} becomes 
\begin{equation}\label{eq:3.6}  
 \frac{\gamma_1(x-e_2)}{\gamma_1(x)} = \frac{\gamma_2(x-e_1)}{\gamma_2(x)}.  
\end{equation} 
Recall that the set $V$ satisfies \eqref{eq:2.3} and $(0,0) \in V$. Clearly, 
multiplying $\gamma_1$ or $\gamma_2$ by a nonzero constant will not change 
the equation, so that we can assume that $\gamma_1(0) =1$ and $\gamma_2(0) 
=1$. For $x_1\ge 0 $ and $x_2 \ge 0$, we define a weight function $W$ by 
$$
W(x) = \gamma_1(x-e_1)\gamma_1(x-2 e_1)\ldots \gamma_1(0, x_2)
       \gamma_2(0, x_2-1)\gamma_2(0, x_2-2) \ldots \gamma_2(0).
$$
Using the equation \eqref{eq:3.6}, we then have
\begin{align*}
\frac{W(x+e_1)}{W(x)} = & \frac{\gamma_1(x)}{\gamma_1(0,x_2)}
  \frac{\gamma_2(1,x_2-1)}{\gamma_2(0,x_2-1)} \ldots
  \frac{\gamma_2(1,1)}{\gamma_2(0,1)} \frac{\gamma_2(1,0)}{\gamma_2(0,0)} \\
  = & \frac{\gamma_1(x)}{\gamma_1(0,x_2)}
  \frac{\gamma_1(0,x_2)}{\gamma_1(0,x_2-1)} \ldots
  \frac{\gamma_1(0,2)}{\gamma_1(0,1)} \frac{\gamma_1(0,1)}{\gamma_1(0,0)} \\
  = & \gamma_1(x). 
\end{align*}
Furthermore, using the equation \eqref{eq:3.7} again, 
\begin{align*}
\frac{W(x+e_2)}{W(x)} = & \frac{\gamma_1(x_1-1, x_2+1)}{\gamma_1(x_1-1,x_2)}
 \frac{\gamma_1(x_1-2,x_2+1)}{\gamma_1(x_1-2,x_2)} \ldots
 \frac{\gamma_1(0,x_2+1)}{\gamma_1(0,x_2)} 
 \frac{\gamma_2(0,x_2)}{\gamma_2(0,0)} \\
 = & \frac{\gamma_2(x)}{\gamma_2(x-e_1)}
 \frac{\gamma_2(x -e_1)}{\gamma_2(x - 2e_1)} \ldots
 \frac{\gamma_2(1,x_2)}{\gamma_2(0,x_2)}\frac{\gamma_2(0,x_2)}{\gamma_2(0,0)}\\
 = & \gamma_2(x). 
\end{align*}
Consequently, the function $W$ defined above is a solution for the equations 
\eqref{eq:3.4}.
\end{proof}

\subsection{Another notion of self-adjointness}
We should like to point out that, just as in the definition of the 
standard difference equation, we can reverse the role of the forward and
the backward differences and give another notion of self-adjointness.

Let $\wt a_{i,j}$ be continuous functions of two variables. We could call
a second order difference operator $\wt \CE$ self-adjoint if 
\begin{equation} \label{eq:3.11}
\wt \CE u =
 \Delta_1\left[\wt a_{1,1}(x)\nabla_1 u + \wt a_{1,2}(x)\nabla_2 u \right]
+\Delta_2\left[\wt a_{2,1}(x)\nabla_1 u + \wt a_{2,2}(x)\nabla_2 u \right].
\end{equation}

For this notion of the self-adjointness, we can also determine the conditions
under which $W \CD$ is self-adjoint. Let us define 
\begin{align*}
  \wt \alpha(x): = & (A_{1,1}(x-e_1)+ B_1(x))(A_{2,2}(x-e_2)+ B_2(x)) - 
     A_{2,1}(x-e_2)A_{1,2}(x-e_1), \\
  \wt r_1(x) := & A_{1,1}(x) + A_{2,1}(x), \quad\hbox{and}\quad
  \wt r_2(x) := A_{1,2}(x) + A_{2,2}(x).
\end{align*}
Following the proof that leads to the result for $\CD$, we can prove the
following:

\begin{proposition}
Let $\CD$ be the difference operator in \eqref{eq:2.5}. There is a function
$W$ such that $W(x) \CD$ is self-adjoint in the sense of \eqref{eq:3.11} if 
and only if $W$ satisfies 
\begin{align} \label{eq:3.12}
\begin{split}
 W B_1 &=\nabla_1(W A_{1,1})(\cdot +e_1) +\nabla_2 (W A_{1,2})(\cdot +e_1)\\
 W B_2 &=\nabla_1(W A_{2,1})(\cdot+e_2) +\nabla_2 (W A_{2,2})(\cdot +e_2).
\end{split}
\end{align}
Furthermore, assume that $\wt \alpha(x)\ne 0$; then there is a function $W$ 
such that $W \CD u$ is self-adjoint if and only if 
\begin{equation*}
\wt \alpha (x -e_1) \wt r_2(x) \wt r_1(x-e_2) = \wt \alpha (x -e_2) 
  \wt r_1(x) \wt r_2(x-e_1). 
\end{equation*}
\end{proposition}

\subsection{Weight function and orthogonality}

Given a weight function $W$ defined on a lattice set $V$, we define 
$$
  \CL f = \sum_{x \in V} f(x) W(x). 
$$
For the difference operator given in $\CD$, we define
\begin{equation}\label{eq:3.7}
  \CA_1 = A_{1,1} \nabla_1 + A_{1,2} \nabla_2 + B_1 \quad \hbox{and}\quad
  \CA_2 = A_{2,1} \nabla_1 + A_{2,2} \nabla_2 + B_2. 
\end{equation} 
From the definition of $\CD$ we have
\begin{equation}\label{eq:3.8}
   \CD u = \CA_1 \Delta_1 u + \CA_2 \Delta_2 u. 
\end{equation}

Let $\partial V$ denote the boundary of the set $V$; that is, $\partial V$ 
is a subset of $V$ such that if $x \in \partial V$ then at least one of the
elements $x \pm + e_i$, $i =1,2$, does not belong to $V$. In order that the
difference equation makes sense on $V$, it needs to satisfy some 
condition on $\partial V$. For example, we can impose the conditions 
\eqref{eq:3.9} below on the coefficients $A_{i,j}$.

\begin{proposition}\label{prop:3.4}
If $W \CD$ is self-adjoint and $A_{i,j}$ satisfy the conditions 
\begin{align} \label{eq:3.9}
\begin{split}
  A_{1,1}(x) & = A_{2,1}(x) = 0, \qquad \hbox{$x \in \partial V$ and $x+e_1 
   \in  V$}\\
  A_{1,2}(x) & = A_{2,2}(x) = 0, \qquad \hbox{$x \in \partial V$ and $x+e_2 
   \in  V$},
\end{split}
\end{align}
then, for any polynomial $g$, 
$$
\CL (\CA_1 g) =0 \quad\hbox{and}\quad \CL (\CA_2 g) =0.
$$ 
\end{proposition}

\begin{proof}
Let $\partial_1 V = \{x \in \partial V: x+e_1 \in V\}$ and 
$\partial_2 V = \{x \in \partial V: x+e_2 \in V\}$. Changing summation 
index and using the fact that $W(x+e_1) = 0$ if $x \in \partial V \setminus 
\partial_1 V$, we obtain
\begin{align*} 
   \CL(A_{1,1} \nabla_1 g)& =\sum_{x \in V} W(x) A_{1,1}(x) \nabla_1 g(x) \\
     & = \sum_{x \in V} W(x) A_{1,1}(x) g(x) - 
            \sum_{x+e_1 \in V} W(x+e_1) A_{1,1}(x+e_1) g(x) \\ 
     & = - \sum_{x \in V} \Delta_1(W A_{1,1})(x) g(x) + 
              \sum_{x \in \partial_1 V} W(x) A_{1,1}(x) g(x). 
\end{align*}
A similar equation holds for $\CL(A_{1,2} \nabla_2 g)$. Hence, by the
assumption on $A_{i,j}$, we conclude that 
\begin{align*} 
\CL(A_{1,1} \nabla_1 g) = & - \sum_{x \in V} \Delta_1(W A_{1,1})(x) g(x), \\
\CL(A_{1,2} \nabla_2 g) = & - \sum_{x \in V} \Delta_2(W A_{1,2})(x) g(x). 
\end{align*}  
Since $W \CD$ is self-adjoint if and only if \eqref{eq:3.1} holds, by 
Proposition \ref{prop:3.3}, it follows that 
\begin{align*}
\CL (\CA_1 g )= & \CL(A_{1,1} \nabla_1 g)+\CL(A_{1,2}\nabla_1 g)+\CL(B_1 g) \\
   = & - \sum_{x \in V} \left[\Delta_1(W A_{1,1})(x)+ \Delta_2(W A_{1,2})(x) 
          - W(x) B_1(x) \right] g(x) = 0.
\end{align*}
The proof that $\CL (\CA_2 g) =0$ is similar. 
\end{proof}

\begin{definition} \label{defn:3.5}
The weight function $W$ is consistent with the difference operator
$\CD$ if, for all $g \in \Pi^2$, 
\begin{enumerate}
\item $\CL (\CA_1 g) =0$ and $\CL (\CA_2 g) =0$; 
\item $\CL (A_{1,2} g(\cdot - e_2)) = \CL (A_{2,1} g(\cdot - e_1))$. 
\end{enumerate}
\end{definition}

In particular, setting $g(x) = x_i$ shows that if $W$ is consistent with 
$\CD$ then
$$
   \CL(B_1) =0 \quad\hbox{and}\quad \CL(B_2) =0.
$$

By Proposition \ref{prop:3.4}, if $W\CD$ is self-adjoint and the 
coefficients of $\CD$ satisfies \eqref{eq:3.9}, then $W$ is consistent
with $\CD$ if the condition (2) in the above definition holds. We need
the additional condition for the following result: 

\begin{theorem}
If $W$ is consistent with the difference operator $\CD$ and the polynomials
$u$ and $v$ satisfy $\CD u = \lambda u$ and $\CD v = \mu v$, then
$$
  (\lambda - \mu) \CL (u v )= 0.
$$
\end{theorem}

\begin{proof}
By the definition of $\CA_1$ and the product formula for $\Delta_j$, 
\begin{align*}
 \CA_1 (f g) & = A_{1,1}[f  \nabla_1 g + g(\cdot - e_1) \nabla_1 f] + 
     A_{1,2} [ f \nabla_2 g + g(\cdot - e_2) \nabla_2 f] + B_1 f g \\
 & =  f \CA_1 g +  g(\cdot - e_1)  A_{1,1} \nabla_1 f + 
         g(\cdot - e_2) A_{1,2} \nabla_2 f. 
\end{align*}
Applying $\CL$ to both sides and using the fact that $W$ is consistent
with $\CD$, we conclude that 
$$
  \CL(f \CA_1 g ) = - \CL( g(\cdot - e_1) A_{1,1} \nabla_1 f + 
         g(\cdot - e_2) A_{1,2} \nabla_2 f ). 
$$
Similarly, working with $\CA_2$, we have
$$
 \CL(f \CA_2 g ) = - \CL( g(\cdot - e_1) \nabla_1 A_{2,2} \nabla_1 f + 
         g(\cdot - e_2) A_{2,2} \nabla_2 f ). 
$$
From these two identities and the fact that  $\Delta_j u(\cdot -e_j) = 
\nabla_j u$, it follows that
\begin{align*}
 \CL (v \CA_1 \Delta_1 u ) & = - \CL (\nabla_1 u A_{1,1} \nabla_1 v) 
    - \CL(\Delta_1 u(\cdot - e_2) A_{1,2} \nabla_2 v),\\
 \CL (v \CA_2 \Delta_2 u ) & =  - \CL(\Delta_2 u(\cdot - e_1) A_{2,1} 
  \nabla_1 v) - \CL (\nabla_2 u A_{2,2} \nabla_2 v).
\end{align*}
Hence, by \eqref{eq:3.8} and the fact that  $\Delta_j u(\cdot -e_j) = 
\nabla_j u$, it follows that
\begin{align*}
(\lambda - \mu )\CL (u v) = & \CL(v \CA_1 \Delta_1 u)+ \CL(v \CA_2 \Delta_2 u)
   - \CL(u \CA_1 \Delta_1 v)-\CL(u \CA_2 \Delta_2 v)\\
= & -\CL(\Delta_1 u(\cdot - e_2)A_{1,2}\nabla_2 v)
  -\CL(\Delta_2 u(\cdot - e_1)A_{2,1}\nabla_1 v)\\
  & +\CL(\Delta_1 v(\cdot - e_2)A_{1,2}\nabla_2 u)
  +\CL(\Delta_2 v(\cdot - e_1)A_{2,1}\nabla_1 u)\\
= & -\CL(A_{1,2} ( \Delta_1 u \Delta_2 v)(\cdot - e_2))
  -\CL(A_{2,1} (\Delta_2 u \Delta_1 v)(\cdot - e_1))\\
  & + \CL(A_{1,2} ( \Delta_1 v \Delta_2 u)(\cdot - e_2))
    + \CL(A_{2,1} (\Delta_2 v \Delta_1 u)(\cdot - e_1)) =0
\end{align*}
by the part (2) in the Definition \ref{defn:3.5} of the consistency.
\end{proof}

\begin{corollary}
If $W \CD$ is self-adjoint and the coefficients of $\CD$ satisfies 
\eqref{eq:3.9} and  
\begin{equation} \label{eq:3.10a}
 \CL (A_{1,2} g(\cdot - e_2)) = \CL (A_{2,1} g(\cdot - e_1)), 
 \qquad g \in \Pi^2,  
\end{equation} 
then the equation $\CD u = \lambda_k u$ has polynomial solutions that are 
orthogonal with respect to $W$ on $V$.
\end{corollary}

\subsection{Condition for consistency}

Let us look at the conditions for the consistency closely. Recall that the 
Type A set $V = \{(x_1,x_2): 0 \le x_1 \le M, 0 \le x_2 \le N\}$ and the 
Type B set $V = \{(x_1,x_2): x_1 \ge 0, x_2 \ge 0,0 \le x_1 +x_2 \le N\}$. 
If $V$ is of Type A or Type B, the condition \eqref{eq:3.9} means 
\begin{align} \label{eq:3.10}
  A_{1,1}(0,x_2) = A_{2,1}(0,x_2) = 0, \quad \hbox{and} \quad
  A_{1,2}(x_1,0) = A_{2,2}(x_1,0) = 0.
\end{align}
Since $A_{1,1}$ and $A_{2,1}$ are coefficients of $\Delta_1 \nabla_1$ and
$\Delta_2 \nabla_1$, both involve $\nabla_1$, and $A_{2,1}$ and $A_{2,2}$
are coefficients of $\Delta_1 \nabla_2$ and $\Delta_2 \nabla_2$, both 
involve $\nabla_2$, the conditions in \eqref{eq:3.10} appear to be reasonable 
assumptions.  

Recall that  $\partial_1 V = \{x \in \partial V: x+e_1 \in V\}$ and 
$\partial_2 V = \{x \in \partial V: x+e_2 \in V\}$. 

\begin{proposition}
Assume that $A_{1,2} - A_{2,1} \in \Pi_1^2$. If $\CD$ is admissible and
\eqref{eq:3.9} holds, then $\partial_1 V$ is a subset of $\{x: x_1 = 0\}$ 
and $\partial_2 V$ is a subset of $\{x: x_2 =0\}$. In particular, 
\eqref{eq:3.9} becomes \eqref{eq:3.10}.
\end{proposition}

\begin{proof}
The condition \eqref{eq:3.9} means that the pair $A_{1,1}(x)$ and 
$A_{2,1}(x)$, and the pair $A_{1,2}(x)$ and $A_{2,2}(x)$, agrees on 
the set $\partial V_1$. Recall that both $A_{1,1}$ and $A_{2,1}$ are 
quadratic polynomials and they are of the form \eqref{eq:2.7}. If 
$\partial V_1$ contains more than 4 points, then the classical Bezout 
theorem on the algebraic curves will imply that $A_{1,1}$ and $A_{2,1}$
has a common factor. Since $A_{1,1}(x) = a x_1^2 + \pi_{1,1}(x)$ and
$A_{2,1}(x) = a x_1 x_2 + \pi_{2,1}(x)$, where $\pi_{i,j}$ are polynomials
of degree $1$,  the common factor has to be a linear polynomial, which 
means that both polynomials are product of linear factors. Thus, they
must be of the form 
$$
 A_{1,1}(x) = a x_1 (x_1 + b_1) \quad \hbox{and} \quad  
 A_{2,1}(x) = a x_1 (x_2 + b_2).  
$$ 
Consequently, $A_{1,1}(0,x_2) =0$ and $A_{2,1}(0,x_1) =0$ become the
only possibility. Thus, $\partial V_1$ is a subset of $\{x: x_1 =0\}$.

The same consideration applies to the pair $A_{1,2}(x)$ and $A_{2,2}(x)$
and shows that $\partial V_2$ is a subset of $\{x: x_2 =0\}$. 
\end{proof}

\begin{corollary} \label{cor:3.10}
Assume that $A_{1,2} - A_{2,1} \in \Pi_1^2$. If $W \CD$ is self-adjoint 
and the weight function $W$ is consistent with $\CD$, then 
\begin{align} \label{eq:3.15}
\begin{split}
& A_{1,1}(x) = a x_1 (c_1 -x_1), \quad  A_{1,2}(x) = a x_2 (c_2 -x_1), \quad 
  B_1(x) = b x + d_1; \\
& A_{2,1}(x) = a x_1 (c_3 -x_2), \quad A_{2,2}(x) = a x_2 (c_4 -x_2), 
  \quad   B_2(x) = b x + d_2.
\end{split}
\end{align}
\end{corollary}

Notice that \eqref{eq:3.10} only imposes conditions on $A_{i,j}$ on part 
of $\partial V$. If $V$ is finite, there should also be conditions on the
set $\partial V \setminus (\partial_1 V \cup \partial_2 V)$. 
The equation \eqref{eq:3.10a} imposes such a condition, which is often 
fulfilled if the weight function $W$ vanishes on the 
set $\partial V \setminus (\partial_1 V \cup \partial_2 V)$. 

The property that $\partial_1 V$ and $\partial_2 V$ are subsets of 
$\{x: x_2 =0\}$ and $\{x: x_1 =0\}$, respectively, is fulfilled by
the Type A and the Type B sets. In fact, they appear to be the only
sets $V$ that satisfy this property in connection with the difference
equation.

\section{Difference equations and orthogonal polynomials}
\setcounter{equation}{0}

\subsection{Preliminary}
The results in the previous section provides us with a method to identify 
difference equations that have orthogonal polynomials as solutions. 
The method consists of the following steps. 

\medskip\noindent
{\bf Step 1}. Among all admissible difference operators that satisfy
\eqref{eq:3.15}, find $\CD$ that satisfies \eqref{eq:3.5}. 

\medskip\noindent
{\bf Step 2}. Using \eqref{eq:3.5a} to identify the weight function that 
satisfies \eqref{eq:3.1}. 

\medskip\noindent
{\bf Step 3}. Verifying that the condition \eqref{eq:3.10a} holds. 

\medskip

We will follow the steps outlined above to do a case by case study. The 
classical discrete orthogonal polynomials of one variable will appear in 
the study, so we list them below (see \cite{NSU} or \cite{KoSw}):

\medskip\noindent
{\bf Hahn polynomials} $u = Q_n(x;\alpha,\beta,N)$:   
$$
 t(N+\beta+1-t) \Delta \nabla u + (N(\alpha+1) - (\alpha+\beta+2)t) \Delta u 
  = - n (n+\alpha+\beta+1) u, 
$$
where $\alpha, \beta \ge 0$ and $N\in \NN$. They are orthogonal with respect 
to the weight function $\binom{\alpha+t}{\alpha}\binom{\beta+N- t}{\beta}$ on 
$V = \{0, 1,\ldots, N\}$.

\medskip\noindent
{\bf Meixner polynomials} $u = M_n(t;\beta,c)$:
$$
  t \Delta \nabla u + [c(t+\beta) - t] \Delta u = n(c-1) u  
$$
where $\beta > 0$ and $0 < c < 1$. They are orthogonal with respect 
to the weight function $\frac{(\beta)_x}{x!} c^x $ on $V= \NN_0$.

\medskip\noindent
{\bf Krawtchouk polynomials} $u = K_n(t;p, N)$:
$$
  t (1-p) \Delta \nabla u + [p(N-t)- t(1-p)] \Delta u = - n u  
$$
where $p > 0$ and $N\in \NN_0$. They are orthogonal with respect to the 
weight function $\binom{N}{t}p^t(1-p)^{N-t}$ on $V= \{1,2,\ldots,N\}$.

\medskip\noindent
{\bf Charlier polynomials} $u = C_n(t;a)$:
$$
  t \Delta \nabla u + (a - t) \Delta u = -n u  
$$
where $a> 0$. They are orthogonal with respect to the weight function
$a^t/t!$ on $V=\NN_0$.

\subsection{Difference equations with orthogonal polynomial solutions}

We will now follow the steps outlined to identify the second order 
difference equations that have orthogonal polynomial as solutions. 
For this we assume that $\CD$ is admissible and the coefficients 
$A_{i,j}$ of $\CD$ satisfy \eqref{eq:3.15}. 

The equation \eqref{eq:3.5} is a nonlinear equation, we solve it 
with the help of a computer algebra system. 

\medskip\noindent
{\bf Case 1: $a \ne 0$}.
We can assume that $a = -1$. Then $A_{i,j}$ and $B_i$ must take the form 
a\begin{align*}
& A_{1,1}(x) = x_1 (c_1 -x_1), \quad  A_{1,2}(x) = x_2 (c_2 -x_1), \quad 
  B_1(x) = b x + d_1; \\
& A_{2,1}(x) =  x_1 (c_3 -x_2), \quad A_{2,2}(x) = x_2 (c_4 -x_2), 
  \quad   B_2(x) = b x + d_2.
\end{align*}
With the $A_{i,j}$ and $B_i$ given as above, solving the equation 
\eqref{eq:3.5} (with the help of a computer algebra system) leads 
to essentially two solutions. 

\medskip
{\bf Case 1 (i)}: One solution is 
$$
d_2 = -(b + c_1 - c_2) c_3, \quad d_1 = c_2 (-b - c_1 + c_2), \quad 
c_4 = c_1 - c_2 + c_3,
$$
so that, by solving \eqref{eq:3.5a}, we get 
$$
W(x) = \frac{\Gamma(-c_2 + x_1)\Gamma(-c_3 + x_2)\Gamma(c_1+c_3-x_1 - x_2)}
    {\Gamma(1 + x_1)\Gamma(1 + x_2)\Gamma(1+b+c_1 -c_2 - x_1 - x_2)}. 
$$ 
It is easy to verify that \eqref{eq:3.10a} 
holds for $V = \{x: x_1 \ge 0, x_2 \ge 0, x_1+x_2 \le N\}$. 

This turns out to be exactly the Hahn polynomials of two variables
discussed in \cite{KM}. Indeed, setting $c_2 = -(\sigma_1 + 1)$, 
$c_3 = -(\sigma_2 + 1)$, $c_1 = N + \sigma_2 + \sigma_3 + 2$ and 
$b = -3 - \sigma_1 - \sigma_2 - \sigma_3$ shows that the difference equation 
is the one given below:

\medskip\noindent
{\bf Hahn Polynomials}. 
\begin{align} \label{eq:4.1}
 & x_1(N-x_1+\sigma_2 + \sigma_3+2) \Delta_1 \nabla_1 u 
    - x_2(x_1+\sigma_1+1) \Delta_1 \nabla_2 u\\
 & \qquad - x_1(x_2+\sigma_2+1) \Delta_2 \nabla_1 u +  
     x_2(N-x_2+\sigma_1 + \sigma_3+2) \Delta_2 \nabla_2 u \notag \\
 & \qquad  
 + \left[(N-x_1)(\sigma_1+1)-x_1(\sigma_2+\sigma_3+2)\right]\Delta_1 u \notag\\
 & \qquad
 + \left[(N-x_2)(\sigma_2+1)-x_2(\sigma_1+\sigma_3+2)\right]\Delta_2 u
 = - n (n + |\sigma| +2) u,  \notag
\end{align}
where $|\sigma| = \sigma_1+\sigma_2+\sigma_3$. The weight function $W$ is 
defined on $V = \{x: x_1 \ge 0, x_2 \ge 0, x_1+x_2 \le N\}$ by 
$$
   W_\sigma(x) = \binom{x_1+\sigma_1}{\sigma_1} \binom{x_2+\sigma_2}{\sigma_2}
       \binom{N-x_1-x_2+\sigma_3}{\sigma_3}.
$$
For each $n$, the equation has solutions $\phi_{l,m}(\cdot;\sigma)$,
$m = n-l$ and $0 \le l \le n$, given in terms of Hahn polynomials by 
\begin{align*} 
 \phi_{l,m}(x;\sigma) = 
  Q_l (x_1;\sigma_1,\sigma_2+\sigma_3+2m+1,N-l) (-N+x)_m
  Q_m(x_2;\sigma_2,\sigma_3,N-x_1).  
\end{align*}

\medskip
{\bf Case 1 (ii)}: The other solution is 
$$
 d_1 = c_2 (-b - c_1 + c_2), \quad d_2 = -(1 + b) c_3, \quad 
   c_4 = 1 + c_3, 
$$
so that, by solving \eqref{eq:3.5a}, we get 
$$
W(x) = \frac{\Gamma(-c_2 + x_1) \Gamma( - 1- b+c_2 + x_2)}
   {\Gamma(1 + x_1)\Gamma(1 + x_2)}. 
$$ 
The consistent condition \eqref{eq:3.10a} with $g(x) = x_i$ shows 
that $b = c_2 + c_3 -1$. Setting $c_2 = - \kappa_1 -1$, $c_3 = -\kappa_2 -1$
and $c_1 = M + \kappa_2+2$ shows that the difference equation is of the 
same form as the equation \eqref{eq:4.1} for the Hahn polynomials of two 
variables with $\sigma_1 = - \kappa_1$, $\sigma_2 = - \kappa_2$ and
$N = - 2 - \kappa_1 - \kappa_2$. Since $M$ has to be a positive integer,
this shows that $\kappa_1$ and $\kappa_2$ need to be negative. However, 
in this case, 
$$
   W(x) =  \frac{\Gamma(1+ \kappa_1 + x_1) \Gamma(1+\kappa_2 + x_2)}
      {\Gamma(1 + x_1)\Gamma(1 + x_2)} 
$$
has singularities at $x_1 = - \kappa_1 -1$ or $x_2 = - \kappa_2 -1$. 
Thus, this is not a proper solution. 

\medskip\noindent
{\it Remark}. The computer algebra system that we used produced several 
solutions for \eqref{eq:3.5}, upon close examination, however, we have 
\eqref{eq:4.1} as essentially the only solution. The other solutions 
either turn out to be special cases of \eqref{eq:4.1} or like the case
(ii) given above. Even though the case (ii) does not lead to a proper 
solution, we have included the case to indicate how the method works.

\medskip\noindent
{\bf Case 2: $a = 0$}.
Assume that $\CD$ is admissible and \eqref{eq:3.10} holds.  Then $A_{i,j}$ 
and $B_i$ must take the form 
\begin{align*}
& A_{1,1}(x) = a_1 x_1, \quad  A_{1,2}(x) = a_2 x_2, \quad 
  B_1(x) = b x + d_1; \\
& A_{2,1}(x) =  a_3 x_1, \quad A_{2,2}(x) =a_4  x_2, 
   \quad   B_2(x) = b x + d_2.
\end{align*}
With these $A_{i,j}$ and $B_i$, solving the equation \eqref{eq:3.5} 
(with the help of a computer algebra system) shows that there are 
several solutions. 

\medskip
{\bf Case 2 (i): $a_2 \ne 0$ and $a_3 \ne 0$}. In this case we get
$$
d_1 = a_2 d_2 / a_3, \quad  a_1 = a_2 - b, \quad a_4 = a_3 - b 
$$
and, by solving \eqref{eq:3.5a}, 
$$
 W(x) = \frac{(-a_2)^x (-a_3)^y (a_2+a_3-b)^{N-x-y}}
    {\Gamma(1-d_2/a_3 - x - y) \Gamma(1 + x)\Gamma(1 + y)}. 
$$ 
It is easy to verify that the condition \eqref{eq:3.10a} holds for 
$V = \{x: x_1 \ge 0, x_2 \ge 0, x_1+x_2 \le N\}$. This turns out 
to be exactly the Krawtchcuk polynomials of two variables on $V$. Indeed,
setting $b = -1$, $a_2 = - p_1$, $a_3 = - p_2$ and $d_2 = a_3(1-N)$
shows that the difference equation is the one given below:

\medskip\noindent
{\bf Krawtchcuk Polynomials}. 
\begin{align} 
 & (1-p_1) x_1 \Delta_1 \nabla_1 u - p_1 x_2 \Delta_1 \nabla_2 u
    - p_2 x_1 \Delta_2 \nabla_1 u + (1-p_2) x_2 \Delta_2 \nabla_2 u \\
 & \qquad  
 + \left[p_1(N-x_1)-(1-p_1) x_1\right]\Delta_1 u 
 + \left[p_2(N-x_2)-(1-p_2) x_2\right]\Delta_2 u \notag \\
& \qquad  = - n u. \notag
\end{align}
The weight function $W$ is defined on
$V = \{x: x_1 \ge 0, x_2 \ge 0, x_1+x_2 \le N\}$ by  
$$
W_p(x) = \frac{N!}{x! y! (N-x-y)!} p_1^{x_1} p_2^{x_2} (1-p_1-p_2)^{N-x_1-x_2}.
$$
For each $n$, the equation has solutions $\phi_{l,m}(\cdot; p,N)$, 
$m = n-l$ and $0 \le l \le n$, given in terms of Krawtchcuk polynomials by
$$
 \phi_{l,m}(x;,p,N) = K_l(x;p_1,N-m) (-N+x)_m K_m(y; p_2/(1-p_1),N-x).
$$ 

These polynomials have appeared in \cite{Tr1} but the difference equation
was not given there.  

\medskip
{\bf Case 2 (ii): $a_2 = 0$ and $a_3 = 0$}. In this case we have 
several solutions, depending on the signs of various coefficients. 
These correspond to the product type orthogonal polynomials, which 
we list below. 

\medskip\noindent
{\bf Product Meixner polynomials:}
\begin{align} 
 &  \frac{1}{c_1-1} x_1 \Delta_1 \nabla_1 u  + 
    \frac{1}{c_2-1} x_2 \Delta_2 \nabla_2 u  \\
 & \qquad  
 + \frac{1}{c_1-1} \left[c_1(x_1 + \beta_1)- x_1\right]\Delta_1 u 
 + \frac{1}{c_2-1} \left[c_2(x_2 + \beta_2)- x_2\right]\Delta_2 u 
   =  n u. \notag
\end{align}
The weight function $W$ is defined on $V = \NN_0^2$ by 
$$
W_{b,\beta}(x) = \frac{(\beta_1)_{x_1}}{x_1!} c_1^{x_1} 
      \frac{(\beta_2)_{x_2}}{x_2!} c_2^{x_2}. 
$$
For each $n$, the equation has solutions $\phi_{l,m}(\cdot; \beta,c)$, 
$m = n-l$ and $0 \le l \le n$, given in terms of Meixner polynomials by
$$
 \phi_n(x;\beta,c) = M_l(x_1;\beta_1,c_1)M_m(x_2;\beta_2,c_2).
$$ 

\medskip\noindent
{\bf Product Meixner-Krawtchcuk polynomials:}
\begin{align} 
 &  \frac{1}{c-1} x_1 \Delta_1 \nabla_1 u  - 
      (1-p) x_2 \Delta_2 \nabla_2 u \\
 & \qquad  
 + \frac{1}{c-1} \left[c(x_1 + \beta)- x_1\right]\Delta_1 u 
 - \left[p (N-x_2) - (1-p) x_2\right]\Delta_2 u 
 =  n u. \notag
\end{align}
The weight function $W$ is defined on $V = \{(i,j): i \in \NN_0, 
0\le j \le N\}$ by 
$$
W_{c,\beta,p}(x) = \frac{(\beta)_{x_1}}{x_1!} c^{x_1} 
      \binom{N}{x_2} p^{x_2} (1-p_2)^{N-x_2}.
$$
For each $n$, the equation has solutions $\phi_{l,m}(\cdot;\beta,c,p,N)$, 
$m = n-l$ and $0 \le l \le n$, given in terms of Meixner and Krawtchcuk 
polynomials by
$$
 \phi_n(x;\beta,c,p,N) = M_l(x_1;\beta,c) K_m(x_2;p, N).
$$ 

\medskip\noindent
{\bf Product Meixner-Charlier polynomials:}
\begin{align} 
 &  \frac{1}{c-1} x_1 \Delta_1 \nabla_1 u  - 
       x_2 \Delta_2 \nabla_2 u  \\ 
  & \qquad+ \frac{1}{c-1} \left[c(x_1 + \beta)- x_1\right]\Delta_1 u 
 -  (a-x_2) \Delta_2 u 
 =  n u. \notag
\end{align}
The weight function $W$ is defined on $V = \NN_0^2$ by 
$$
W_{c,\beta,a}(x) = \frac{(\beta)_{x_1}}{x_1!} c^{x_1} 
      \frac{(a)_{x_2}}{x_2!}.
$$
For each $n$, the equation has solutions $\phi_{l,m}(\cdot; \beta,c,a)$, 
$m = n-l$ and $0 \le l \le n$, given in terms of Meixner and Charlier 
polynomials by
$$
 \phi_n(x;\beta,c,a) = M_l(x_1;\beta,c) C_m(x_2;a).
$$ 

\medskip\noindent
{\bf Product Krawtchcuk polynomials:}
\begin{align} 
 &  (1-p_1) x_1 \Delta_1 \nabla_1 u  + 
    (1-p_2) x_2 \Delta_2 \nabla_2 u  \\
 &   
 + \left[p_1(N_1-x_1)- (1-p_1)x_1 \right]\Delta_1 u 
 + \left[p_2(N_2-x_2)- (1-p_2)x_2 \right]\Delta_2 u \notag = - n u. \notag
\end{align}
The weight function $W$ is defined on $V=\{(i,j): 0\le i \le N_1,
0\le j\le N_2\}$ by 
$$
W_{p}(x) = \binom{N_1}{x_1} p_1^{x_1}(1-p_1)^{N_1-x_1} 
      \binom{N_2}{x_2} p_2^{x_2}(1-p_2)^{N_2-x_2}.
$$
For each $n$, the equation has solutions $\phi_{l,m}(\cdot;p,N)$, 
$m = n-l$ and $0 \le l \le n$, given in terms of Krawtchcuk polynomials by
$$
 \phi_n(x;p,N) = K_l(x_1;p_1,N_1)K_m(x_2;p_2,N_2).
$$ 

\medskip\noindent
{\bf Product Krawtcheuk-Charlier polynomials:}
\begin{align} 
 &  (1-p) x_1 \Delta_1 \nabla_1 u  + 
     x_2 \Delta_2 \nabla_2 u  \\
 & \qquad  
 + \left[p(N-x_1)- (1-p)x_1 \right]\Delta_1 u 
 + (a-x_2) \Delta_2 u  = - n u. \notag
\end{align}
The weight function $W$ is defined on $V=\{(i,j): 0\le i \le N, j \in \NN_0\}$
by 
$$
W_{p,a}(x) = \binom{N,x_1} p^{x_1}(1-p)^{N-x_1} \frac{(a)_{x_2}}{x_2!}
$$
For each $n$, the equation has solutions $\phi_{l,m}(\cdot;a,p,N)$, 
$m = n-l$ and
$0 \le l \le n$, given in terms of Krawtchcuk and Charlier polynomials by
$$
 \phi_n(x;a,p,N) = K_l(x_1;p, N) C_m(x_2;a).
$$ 

\medskip\noindent
{\bf Product Charlier polynomials:}
\begin{align} 
    x_1 \Delta_1 \nabla_1 u  + 
     x_2 \Delta_2 \nabla_2 u 
 + (a_1-x_1)\Delta_1 u 
 + (a_2-x_2)\Delta_2 u  = - n u. 
\end{align}
The weight function $W$ is defined on $\NN_0^2$ by 
$$
W_a(x) = \frac{a_1^{x_1} }{x_1!}\frac{a_2^{x_2} }{x_2!}
$$
For each $n$, the equation has solutions $\phi_{l,m}(\cdot;a)$, $m = n-l$ and
$0 \le l \le n$, given in terms of Charlier polynomials by
$$
 \phi_n(x;a) = C_l(x_1;a_1)C_m(x_2;a_2).
$$ 

\medskip 
{\bf Case 2 (iii), other possibilities}. There is no solution if one of the 
$a_2$ or $a_3$ is zero and the other one is not. However, if we set $a_4 =0$, 
then solving the equation \eqref{eq:3.5} leads to 
$$
  d_1 = -(a_1+b)d_2/b, \quad a_2 = a_1+b, \quad a_3 =b.
$$ 
Solving the equation \eqref{eq:3.5a} leads to the weight function
$$
  W(x) = b^y (a_1+b)^{N-y} \Gamma(d_2/b + x+y) / (x! y!). 
$$
Since $W$ does not have enough decay, the set $V$ needs to be bounded 
to keep the linear functional $\CL$ defined by $W$ well-defined on 
polynomials. However, we do not see a way to choose $V$ such that 
$W$ vanishes on the boundary $\partial V \setminus
(\partial_1 V \cap\partial_1 V)$. In other words, we do not see a way
to choose $V$ such that the condition \eqref{eq:3.10a} can be satisfied. 

\medskip

\subsection{Final Remarks}
  
Under the assumption that $A_{i,j}$ are of the form \eqref{eq:3.15}, we 
found eight difference equations that have orthogonal polynomials as 
solutions. These are the Hahn and the Krawtchcuk polynomials on Type A 
set $V$, and product type polynomials on Type B set $V$: 
product Meixner, product Krawtchcuk, and product Charlier polynomials, 
as well as product Meixner-Krawtchcuk, Meixner-Charlier, and
Krawtchcuk-Chalier polynomials. With the restrictions that we put
on the difference equations (that is, the conditions in Corollary 
3.10), these appear to be the only possible solutions. 

Let us mention that another attempt of characterizing orthogonal
polynomials of two variables satisfying second order difference 
equations,  based on a matrix approach, is currently being undertaken 
by a group in University of Granada (\cite{MFPP}). 

As pointed out in the introduction, if $u$ satisfies a difference equation
then the function $u_h(x):= u(h_1 x_1, h_2 x_2)$ satisfies a similar 
difference equation which becomes in limit, as $h \to 0$, the second order 
partial differential equation. The second order partial differential
equations that have orthogonal polynomials as solutions are characterized
in \cite{KS}. It is known that there are nine such equations, whose 
solutions correspond to product type orthogonal polynomials of two variables, 
as well as two other type, one is the Jacobi type orthogonal polynomials
on the triangle $T=\{(x,y): x \ge 0, y \ge 0, 0 \le x+y \le 1\}$ and the
other is the orthogonal polynomials on the disk $B=\{(x,y): x^2+y^2 \le 1\}$.
The discrete Hahn polynomials of two variables on Type A set appear to
become the Jacobi type polynomials on the triangle. However, as far as 
we can tell, there is no discrete analog for the orthogonal polynomials 
on the disk. 

A natural question is if there are other difference equations of the 
form \eqref{eq:1.3} that have orthogonal polynomials as solutions. 
Notice that our results are obtained under various assumptions on
$\CD$. For example, we assume that the quadratic parts of $A_{1,2}$
and $A_{2,1}$ are equal, that is, $A_{1,2} - A_{2,1} \in \Pi_1^2$. 
Furthermore, we establish the orthogonality by requiring that $\CD$
is self-adjoint and $W$ is consistent with $\CD$. These assumptions
are sufficient but may not be necessary. Still, assuming that
$A_{1,2} - A_{2,1} \in \Pi_1^2$, it seems to us that likely no other 
nontrivial difference equations can have orthogonal polynomials 
as solutions.

\enddocument